\documentclass[fleqn]{mat01}
\usepackage{times,mathtimy,amssymb}
\begin{document}

\font\xx=msam5 at 10pt
\def\ab{\mbox{\xx{\char'03}}}

\setcounter{page}{87}
\firstpage{87}


\def\defi{\trivlist\item[\hskip\labelsep{\bf DEFINITION.}]}
\def\remark{\trivlist\item[\hskip\labelsep{\it Remark.}]}
\def\remarks{\trivlist\item[\hskip\labelsep{\it Remarks}]}
\def\lema{\trivlist\item[\hskip\labelsep{{\it Lemma} A.}]}
\def\noot{\trivlist\item[\hskip\labelsep{{\it Note.}}]}
\def\thet{\trivlist\item[\hskip\labelsep{{\bf Theorem AR.}}]}
\newtheorem{theo}{Theorem}[section]
\renewcommand\thetheo{\arabic{section}.\arabic{theo}}
\newtheorem{theor}{\bf Theorem}
\newtheorem{lem}{Lemma}
\newtheorem{coro}{\rm COROLLARY}
\newtheorem{pro}{\it Proof of Theorem}
\def\ds{\displaystyle}
\def\cl #1{\overline #1}

\title{Some approximation theorems}

\markboth{N V Rao}{Some approximation theorems}

\author{N V RAO}

\address{Department of Mathematics, The University of Toledo, Toledo, Ohio 43606, USA\\
\noindent E-mail: rnagise@math.utoledo.edu}

\volume{113}

\mon{February}

\parts{1}

\Date{MS received 25 February 2002; revised 24 December 2002}

\begin{abstract}
The general theme of this note is illustrated by the following theorem:\vspace{.6pc}

\noindent {\bf Theorem 1.}\ {\it Suppose $K$ is a compact set in the complex plane
and $0$ belongs to the boundary $\partial K$. Let ${\cal A}(K)$ denote
the space of all functions $f$ on $K$ such that $f$ is holomorphic in a
neighborhood of $K$ and $f(0)=0$. Also for any given positive integer
$m$, let ${\cal A}(m,K)$ denote the space of all $f$ such that $f$ is
holomorphic in a neighborhood of $K$ and
$f(0)=f^{\prime}(0)=\cdots=f^{(m)}(0)=0$. Then ${\cal A}(m,K)$ is dense
in ${\cal A}(K)$ under the supremum norm on $K$ provided that there
exists a sector $W=\{r\hbox{\rm e}^{i\theta};\;0\leq
r\leq\delta,\alpha\leq\theta\leq\beta\}$ such that $W\cap K=\{0\}$.
 (This is the well-known Poincare's external cone condition).}\vspace{.5pc} We present
various generalizations of this result in the context of higher
dimensions replacing holomorphic with harmonic.
\end{abstract}

\maketitle

\section{Introduction}

Axler and Ramey (personal communication) have obtained the following
interesting result: Let $L^2(S^n)$ denote the usual Lebesgue space on
the unit sphere $S^n$ with respect to the surface area measure on $S^n$;
$x_0$ be a fixed point in $R^n$. Let ${\cal P}(x_0,m)$ denote the space
of all harmonic polynomials which vanish at $x_0$ together with all
their derivatives of order less than or equal to $m$, and $m>0$. Then\renewcommand\thefootnote{}
\footnote{Dedicated to Prof. Ashoke Roy on his 62nd birthday.} \vspace{.5pc}

\begin{thet}
{\it ${\cal P}(x_0,m)$ is dense in $L^2(S^n)$ if
and only if $|x_0|\geq 1$.}
\end{thet}

\noindent They also posed the following questions:
\begin{enumerate}
\renewcommand{\labelenumi}{(\arabic{enumi})}
\item  Does the above result remain valid if $L^2(S^n)$ is replaced by
any $L^p(S^n)$ with $p>2$?
\item  Could $S^n$ be replaced by more general surfaces?\vspace{-.3cm}
\end{enumerate}

We shall show here that the answer to the 1st question is yes. Let
${\cal C}(x_0,S^n)$ denote the space of all continuous functions on
$S^n$ that vanish at $x_0$ and the space of all continuous functions on
$S^n$.

\setcounter{theor}{1}
\begin{theor}[\!]
For any positive integer $m${\rm ,} ${\cal P}(x_0,m)$ is dense in ${\cal
C}(x_0,S^n)$ with the sup norm if and only if $|x_0|\geq 1$. When $x_0$
is not on the sphere{\rm ,} then ${\cal C}(x_0,S^n)$ is the same as ${\cal
C}(S^n)$. 
\end{theor}\pagebreak

\begin{remark}
We shall not prove that the density fails when $|x_0|<1$ since it is
rather obvious and we shall not explicitly deal with the case when
$|x_0|>1$, because the proof for $|x_0|=1$ can be imitated without any
problems.
\end{remark}

We will derive Theorem~2 as a corollary of a more general result for
which we need to introduce some more notation. Let $K$ be any compact
set in $R^n$, $\partial K$ its boundary. We define a notion called 
ECC. (This is the well-known Poincare's external cone condition.) We
say that {\it K satisfies} ECC {\it at a point} $x_0$ {\it if there exists a
closed solid truncated cone $W$ with vertex at $x_0$ such that $W\cap
K=\{x_0\}$.} 

It is clear that to satisfy ECC at $x_0$, $x_0$ must be on the
boundary of $K$ and also that the set of points where $K$ satisfies 
ECC is dense in the boundary of $K$. In order to see this, take any
point $\xi$ on the boundary of $K$ and a ball of radius $r$ with center
at $\xi$, where $r$ is arbitrary and positive. There must exist a point
$\eta$ outside $K$ such that $|\eta-\xi|<r/2$ for otherwise $\xi$ would
be an interior point of $K$. Now choose a nearest point to $\eta$
in $K$, say $\lambda$. Clearly $|\lambda-\eta|=\delta<r/2$ and the ball
of radius $\delta$ with center at $\eta$ is entirely contained in the
ball of radius $r$ with center at $\xi$. Now $\lambda$ must belong to
the boundary of $K$, must lie within a distance of $r$ from $\xi$ and
satisfies {ECC} for $K$. 

Let ${\cal H}(x_0,K)$ denote the space of all functions $f$ on $K$ such
that $f$ vanishes at $x_0$ and is the restriction to $K$ of a function
harmonic in a neighborhood of $K$. Let ${\cal H}(m,x_0,K)$ denote the
space of all functions $f$ on $K$ such that $f$ is the restriction to
$K$ of a function, harmonic in a neighborhood of $K$ and it, together
with all its derivatives of order $\leq m$ vanish at $x_0$. 

We shall assume the following well-known result:

\begin{lema}
{\it Let $K$ be any closed ball in $R^n$ and $x_0$ belong to $K$.
Then ${\cal P}(x_0,m)$ is dense in ${\cal H}(m,x_0,K)$.
}
\end{lema}\vspace{.5pc}

Also we need

\begin{theor}[\!]
Assume $K$ satisfies  {\rm ECC} at $x_0$. Then for any positive integer
$m$, $\overline{{\cal H}(m,x_0,K)}$ $\supset {\cal H}(x_0,K)$ with the sup
norm.
\end{theor}

We shall supply a proof of this later.
\setcounter{pro}{1}
\begin{pro}
{\rm Let $K$ be the closed unit ball in $R^n$ and
$x_0$ belong to $S^n=\partial K$. Certainly $K$ satisfies {ECC} at
$x_0$. Let $f$ belong to ${\cal C}(x_0,S^n)$. Let $\varepsilon$ be any
positive number. It is well-known that there exists a harmonic
polynomial $P$ such that
\begin{equation*}
|f(x)-P(x)|<\varepsilon\hbox{\quad on\quad} S^n.
\end{equation*}
 Let $h(x)=P(x)-P(x_0)$. Then $|f(x)-h(x)|\leq
|f(x)-P(x)|+|P(x_0)|<2\varepsilon$ on $S^n$ and also $h(x)$ belongs to
${\cal H}(x_0,K)$. But by Theorem~3, there exists a $g$ in ${\cal
H}(m,x_0,K)$ such that $|h(x)-g(x)|<\varepsilon$ and so
$|f(x)-g(x)|<3\varepsilon$. This proves Theorem 2 in view of Lemma A.
}
\hfill \ab
\end{pro}

Since ${\cal C}(x_0,S^n)$ is dense in all $L^p(S^n)$ for $0<p<\infty$, from
Theorem 2, we have

\setcounter{coro}{3}
\begin{coro}$\left.\right.$\vspace{.3pc}

\noindent For any $p,\;0<p<\infty${\rm ;} for any positive integer
$m${\rm ,} and any point $x_0$ on $S^n${\rm ,} the space ${\cal P}(m,x_0)$ of harmonic
polynomials that vanish together with all their derivatives of order less than
or equal to $m$ is dense in $L^p(S^n)$.
\end{coro}

\begin{pro}
{\rm Let $G(x)=\ln|x|$, if $n=2$ and $ |x|^{2-n}$,
if $n>2$. We may assume without loss of generality that $x_0=0$,
$W=\{z;|z|\leq\rho,{z/|z|}\in\hbox{\rm a spherical cap}\ D\}$,
and $W\cap K=\{0\}$.

Fix a $z$ outside $K$. Then $G(x-z)$ is harmonic as a function of $x$ in
a neighborhood of $K$ and in a neighborhood of the origin can be
expanded in an absolutely convergent power series
\begin{equation*}
G(x-z)=\sum a_{\alpha}(z) x^{\alpha}
\end{equation*}
where $\alpha$ is a multi-index $(\alpha_1,\alpha_2,\ldots,\alpha_n)$
and indices are allowed to run through all non-negative integers. Let
$|\alpha|$ denote $\alpha_1+\alpha_2+\cdots+\alpha_n$. Further we notice
that for any fixed non-negative integer $k$, the polynomial
$\sum_{|\alpha|=k}a_{\alpha}(z) x^{\alpha}$ is harmonic in $x$ and for
any fixed $\alpha$, $a_{\alpha}(z)$ is harmonic in $z$ except at the
origin and if we set $z=|z|\omega$ where $\omega$ varies on the unit
sphere,
\begin{equation}
a_{\alpha}(z)=|z|^{2-n-|\alpha|}a_{\alpha}(\omega).
\end{equation}
We note that $a_{\alpha}(\omega)$ is real-analytic on the unit sphere. Now let
\begin{equation*}
G(m,x,z)=G(x-z)-\sum_{|\alpha|\leq m}a_{\alpha}(z)x^{\alpha}.
\end{equation*}
Clearly for any fixed $z\neq o$, $G(m,x,z)$ is harmonic in a
neighborhood of $K$ and vanishes together with its derivatives of order
less than or equal to $m$. Then if $\mu$ is any finite Borel measure on
$\partial K$ orthogonal to ${\cal H}(m,0,K)$, it follows that
\begin{equation*}
\int G(x-z)\,\hbox{d}\mu(x)=\sum_{|\alpha|\leq m}a_{\alpha}(z)\int
x^{\alpha}\ \hbox{d}\,\mu(x)\quad \begin{array}{l}
\hbox{\rm for all}\ z\ \hbox{\rm in}\\ \hbox{\rm the complement of}\
K.\end{array}\tag{0}
\end{equation*}
Let $b_{\alpha}$ denote $\int x^{\alpha}\,\hbox{d}\mu(x)$ and
$p(k,z)$ denote $\sum_{|\alpha|=k}a_{\alpha}(z)b_{\alpha}$. By (1)
it follows that $p(k,z)=|z|^{2-n-k}p(k,\omega)$ and $p(k,\omega)$ is
real-analytic on the unit sphere. We claim that
\begin{equation}
p(k,z)\equiv 0 \quad \hbox{\rm for all}\ \ k,1\leq k\leq m.
\end{equation}
Suppose not.
Then there would exist a positive integer $l$ such that $p(j,z)\equiv
0\hbox{\rm\ for\ }j>l$ $\hbox{\rm\ and\ }p(l,z)\neq 0$.
Because $p(l,z)$ is homogeneous and is real-analytic, the set of its
zeroes on the unit sphere would be a closed set without any interior.
Hence there would exist sub-cone $V$ of $W$ and a positive number
$\delta$ such that
\begin{equation}
|p(l,z)|\geq\delta|z|^{2-n-l}\hbox{\rm\quad for
all\ \ } z\in V
\end{equation}
and further by choosing a sufficiently small
$\beta<\rho$, we have
\begin{equation}
\left|\sigma(z)=\sum\limits_{0\leq k\leq l} p(k,z)\right|\geq \frac{\delta}{2}|z|^{2-n-l} \hbox{\rm\quad
on\ }U=V\cap\{|z|\leq\beta\}.
\end{equation}
Choose a hyper-plane section
$S$ of $U$ through the origin and integrate $\sigma(z)$ on $S$ with
respect to the surface measure on it. Since $\sigma(z)$ stays away from
$0$ on $U$, it has the same sign everywhere and so from (4) it follows that
\begin{equation}
\left\vert\int_{S}\sigma(z)\,\hbox{d}z\right\vert\geq\int|\sigma(z)|\,\hbox{d}z\geq \frac{\delta}{2}
\int_{S}|z|^{2-n-l}\,\hbox{d}z.
\end{equation}
But the last integral is infinite for
$l>0$. But on the other hand $\int_{S}|G(x-z)|\,\hbox{d}z$ is uniformly bounded
and so $\int_{S}|\int G(x-z) \,\hbox{d}\mu(x)|\,\hbox{d}z$ is finite. This, (5), and
(0) lead to a contradiction establishing (2). Hence
\begin{equation}
\int G(x-z)\,\hbox{d}\mu(x)=b_0G(z)\hbox{\rm\quad for every\ }z\hbox{\rm \ outside}\ K.
\end{equation}
If $\nu=\mu-b_0\delta_0$ where $\delta_0$ is the
Dirac measure at the origin, (6) can be restated as
\begin{equation}
\int
G(x-z)\,\hbox{d}\nu(x)=0\hbox{\rm\quad for all\ }z\hbox{\rm\
outside\ }K.
\end{equation}
(7) implies $\nu$ is orthogonal to any
function $f$ which is the restriction to $K$ of a function harmonic in a
neighborhood of $K$. This is a rather standard Runge argument and we omit
the proof. Hence for any $f$ in ${\cal H}(0,K)$, $\int
f(x)\,\hbox{d}\nu(x)=\int f(x)\,\hbox{d}\mu(x)-b_0f(0)=0$ and so $\int
f(x)\,\hbox{d}\mu(x)=0$. Now by Hahn--Banach, we have Theorem 3.
}
\hfill \ab
\end{pro}

\setcounter{pro}{0}
\begin{pro}
{\rm 
Fix a $z$ outside $K$ and write the Taylor
formula of order $m$ for the Cauchy kernel: 
\begin{equation}
\frac{1}{x-z}=-\sum\limits_{0\leq
k\leq m} \frac{x^n}{z^{n+1}}\,+ \frac{x^{m+1}}{z^{m+1}(x-z)}.
\end{equation}
Let $\mu$ be any finite Borel measure on $\partial K$ such that
\begin{equation}
\int f(x)\,\hbox{d}\mu(x)=0\quad \hbox{\rm for any\ } f\in {\cal A}(m,K).
\end{equation}
So $\int{x^{m+1}/{z^{m+1}(x-z)}}\,\hbox{d}\mu(x)=0$ and consequently
\begin{equation}
\int \frac{1}{(x-z)}=-\sum_{0\leq k\leq m}\frac{\int x^n\,\hbox{d}\mu(x)}
{z^{n+1}}.
\end{equation}
 Let ${a}_k=\int x^k\,\hbox{d}\mu(x),0\leq k\leq m$. Arguing
as in the proof of Theorem 3, we find that $a_k=0,1\leq k\leq m$ and
$\mu-a_0\delta_0$ is orthogonal to all functions holomorphic in a
neighborhood of $K$ and so to ${\cal A}(0,K)$. But $\delta_0$ is
orthogonal to ${\cal A}(K)$ and hence follows the theorem.
}
\hfill \ab
\end{pro}

\section{Conclusion}
Several problems remain. One of them is whether
{ECC} is really necessary. Another one is what is the capacity of
the set of points where the conclusion of either Theorem 3 or Theorem 1
holds in analogy with the set of regular points for the Dirichlet
problem? Lastly, what would be an analogue of this Theorem 1 in the
context of several complex variables?

\end{document}